\documentclass{amsart}
\usepackage[fontsize=12pt]{fontsize}
\usepackage{blindtext}
\usepackage{amsmath,amssymb,amsthm,hyperref}
\usepackage{geometry}
\usepackage{bm}
\usepackage{graphicx} 
\usepackage{todonotes}

\def\IND{\setbox0=\hbox{$x$}\kern\wd0\hbox to 0pt{\hss$\mid$\hss}
\lower.9\ht0\hbox to 0pt{\hss$\smile$\hss}\kern\wd0}
\def\NotIND{\setbox0=\hbox{$x$}\kern\wd0\hbox to 0pt{\mathchardef\nn=12854\hss$\nn$\kern1.4\wd0\hss}\hbox to 0pt{\hss$\mid$\hss}\lower.9\ht0 \hbox to 0pt{\hss$\smile$\hss}\kern\wd0}

\def\ind{\mathop{\mathpalette\IND{}}}

\def\indi#1{\mathop{\mathpalette\IND{}^{\!\!\!\!\rlap{$\scriptstyle{#1}$}\,\,\,\,}}}

\theoremstyle{definition}
\newtheorem*{thma}{Main Theorem}
\newtheorem{thm}{Theorem}[section]
\newtheorem{lemma}[thm]{Lemma}
\newtheorem{fact}[thm]{Fact}
\newtheorem{prop}[thm]{Proposition}
\newtheorem{rmk}[thm]{Remark}
\newtheorem{defi}[thm]{Definition}
\newtheorem{coro}[thm]{Corollary}
\newtheorem*{nota}{Notation}
\newtheorem*{rmk*}{Remark}
\newtheorem*{exmp}{Example}

\newcommand{\A}{\mathbb{A}}
\newcommand{\eq}{^\mathrm{eq}}
\newcommand{\alg}{^\text{alg}}

\newcommand{\tp}[1]{\text{tp}(#1)}
\newcommand{\tP}[1]{\text{tp}_{L_P}(#1)}
\newcommand{\acl}[1]{\text{acl}(#1)}

\newcommand{\aclP}[1]{\text{acl}_P(#1)}

\newcommand{\aeq}[1]{\text{acl}^\mathrm{eq}(#1)}
\newcommand{\aeqs}[1]{\text{acl}^\text{eq}_{P,\sigma}(#1)}
\newcommand{\aeqP}[1]{\text{acl}^\text{eq}_{P}(#1)}
\newcommand{\deqP}[1]{\text{dcl}^\text{eq}_{P}(#1)}

\newcommand{\Cb}[1]{\text{Cb}(#1)}
\newcommand{\code}[1]{\lceil #1 \rceil}

\title{Geometric Classification of Imaginaries in Lovely Pairs of Pseudofinite Fields}
\author{Zixuan ZHU}
\thanks{The author has been supported by the China Scholarship Council, the Deutsche Forschungsgemeinschaft (DFG, German Research Foundation) under the Excellence Strategy EXC 2044–390685587, Mathematics Münster: Dynamics–Geometry–Structure, and via the DFG project 495759320 Modelltheorie bewerteter Körper mit Endomorphismus.}
\address{Department for Mathematical Logic and Foundational Research, University of Münster}
\email{zixuan.zhu@uni-muenster.de}

\begin{document}

	\begin{abstract}
In \cite{PILLAY200713}, Pillay showed that every imaginary in a proper pair of algebraically closed fields is interalgebraic with an imaginary of a certain algebraic form. In this paper, we establish an analogue for lovely pairs of pseudofinite fields. After naming countably many constants, we prove every imaginary in lovely pairs of pseudofinite fields is interalgebraic with an imaginary of the same algebraic form as \cite{PILLAY200713} in the corresponding pairs of finite extensions of pseudofinite fields.
\end{abstract}
    \maketitle
	\section*{Introduction}
    In \cite{MR704080}, Poizat showed that, for a complete stable $L$-theory $T$, there is a complete theory $T_P$ in the expanded language $L_P:=L\sqcup \{P\}$ whose models form a certain class of rich elementary pairs of models of $T$, called \emph{beautiful pairs}. Later, in \cite{BENYAACOV2003235}, Ben Yaacov, Pillay, and Vassiliev generalized this construction to complete simple theories,  introducing the corresponding notion of lovely pairs. The complete $L_P$-theory $T_P$ of all lovely pairs is only model-theoretically meaningful if it is \emph{first order}, i.e., if sufficiently saturated models of $T_P$ are lovely pairs..

Even when $T$ eliminates imaginaries, the pair theory $T_P$ may have new imaginaries that cannot be eliminated in the original sorts. In \cite{PILLAY200713}, Pillay described the imaginaries in beautiful pairs (i.e., proper pairs) of algebraically closed fields. Namely, he showed that every imaginary is interalgebraic with an imaginary of a certain algebraic form, which we call \emph{Pillay form} in this paper.

In this paper, we consider the imaginaries in lovely pairs of pseudofinite fields. Recall that a pseudofinite field is an infinite field satisfying every first-order sentence that holds in all finite fields. The theory Psf of pseudofinite fields is a simple theory, and for every completion T of Psf, $T_P$ is first-order.

Our approach is to realize lovely pairs of pseudofinite fields as fixed fields of lovely pairs of models of ACFA. The latter can also be viewed as beautiful pairs of algebraically closed fields equipped with a generic automorphism. This allows us to use Pillay's results on imaginaries in beautiful pairs of algebraically closed fields to study imaginaries in lovely pairs of pseudofinite fields. Our main result is the following.
	\begin{thma}
		After adding countably many constants to the language, lovely pairs of pseudofinite fields admit geometric elimination of imaginaries in Pillay sorts in the corresponding pairs of finite extensions.
	\end{thma}

    We also establish a reduction principle (see Proposition~\ref{P:reductprncpl}) for studying imaginaries in lovely pairs of pseudofinite fields: questions concerning algebraic closure and non-forking independence of imaginaries can be reduced to the corresponding questions in pairs of algebraically closed fields, where these notions have been characterized in \cite{zhu2026}. Our approach also yields a geometric classification of imaginaries in lovely pairs of ACFA, in terms of Pillay forms in the underlying pairs of algebraically closed fields (see Proposition~\ref{P:GEIinK}).
    
    The paper is organized in three sections: In the preliminaries, we recall the necessary background on lovely pairs, stable theories with a generic automorphism, pseudofinite fields, and Pillay's geometric description of imaginaries in pairs of algebraically closed fields. In the second section, we realize lovely pairs of pseudofinite fields as fixed fields of lovely pairs of models of ACFA and obtain representations of their imaginaries by Pillay forms in pairs of finite extensions, with interalgebraicity computed in the ambient structure. We then interpret these finite extensions in the original pair after naming suitable parameters. In the last section, we describe the induced structure on the fixed-field pair and transfer this interalgebraicity to the original pair, completing the proof of the Main Theorem. We also establish the reduction principle for algebraic closure and non-forking independence.
	
	\section{Preliminaries}
	In this section, we recall some concepts and results that are needed in the paper. We start with the notion of lovely pairs.
	\subsection*{Lovely pairs}
	In this subsection, we fix a complete simple $L$-theory $T$ with QE. To study the theory of elementary pairs of models of $T$, we expand the language by adding a unary predicate symbol $P$ naming an elementary submodel. Denote by $L_P := L \sqcup \{P\}$ the expanded language of pairs.
	\begin{defi}
		We say an $L_P$-structure $(M,P)$ is a \textbf{lovely pair} of models of $T$ if $M\succ P(M)\models T$ and if the following properties are satisfied:
		\begin{enumerate}
			\item  For any small\footnote{Here, small means the cardinality of $A$ is no larger than $|T|$.} parameter set $A \subset M$ and any $L$-type $p \in S(A)$, there exists a non-forking extension of $p$ over $A \cup P(M)$ that is realized in $M$. 
			\item For any small parameter set $A \subset M$ and any $L$-type $p \in S(A)$, if $p$ does not fork over $A \cap P(M)$, then $p$ is realized in $P(M)$. 
		\end{enumerate}
	\end{defi}
	\begin{rmk*}
		The first property is called the \emph{heir property}, and the second the \emph{coheir property}.
	\end{rmk*}
	When $T$ is stable (and countable), lovely pairs of models of $T$ coincide with the \emph{beautiful pairs} (of models of $T$) introduced by Poizat in \cite{MR704080}.
	
	In \cite[Corollary 3.9]{BENYAACOV2003235}, it is shown that if $T$ is simple, then lovely pairs always exist and their common theory is complete. We denote this theory by $T_P$. Moreover, if every (equivalently, some) $|T|^+$-saturated model of $T_P$ is a lovely pair of models of $T$, then $T_P$ is said to be \textbf{first-order}. In fact, for every theory $T$ considered in this article, $T_P$ will be first-order.
	In order to be able to state a quantifier reduction result, we first introduce some notation for $T_P$.
	
	For $A\subset(M,P)\models T_P$, we denote by $\acl{A}$ (resp. $\aclP{A}$) the algebraic closure of $A$ in the sense of $T$ (resp. in the sense of $T_P$), and by $\ind$ (resp. $\indi{P}\,$) non-forking independence in the sense of $T$ (resp. $T_P$). Let $P(A):=A\cap P(M)$. By abuse of notation, we also write $P$ for $P(M)$. We say a small parameter set $A$ is \textbf{$P$-independent} if $A\ind_{P(A)}P$. 
	
	After we add a new predicate to the language, the structure is enriched. However, for a $P$-independent set $A$, it suffices to consider the $L$-structure of $A$, as shown by the following proposition that the $L$-structure of a $P$-independent set $A$ essentially determines its $L_P$-structure.
	\begin{fact}[{\cite[3.8]{BENYAACOV2003235}}, {\cite[2.5]{Imgnrinbtfprs}} and \cite{MR704080}]\label{F:TP}
		Assume $T_P$ is first-order. Let $(M,P)\models T_P$. Let $A,B\subset(M,P)$ be small such that $A\equiv^L B$ and both $A$ and $B$ are $P$-independent. Then, we have:
		\begin{enumerate}
			\item If $A$ and $B$ have the same quantifier-free $L_P$-type, then $A$ and $B$ have the same $L_P$-type.
			\item The algebraic closure in the sense of $T_P$ coincides with that in $T$, i.e., $\aclP{A}=\acl{A}$.
			\item If $T$ is stable, then so is $T_P$. Actually, $T_P$ and $T$ have the same stability spectrum.
            \item If $T$ is simple (resp. supersimple), then $T_P$ is also simple (resp. supersimple). 
		\end{enumerate}
		
	\end{fact}
	\subsection*{Stable theories with a generic automorphism}
	Let $T$ be a complete stable $L$-theory. Similarly, to study the theory $T$ with a generic automorphism, we add a unary function symbol $\sigma$ to the language. The new language is denoted by $L_\sigma=L\sqcup\{\sigma\}$.
	
	Now consider the theory $T_\sigma$ whose models are structures $(M,\sigma)$ where $M$ is a model of $T$ and $\sigma$ is an $L$-automorphism of $M$. If $T_\sigma$ admits a model companion, it will be denoted by $T$A. In this case, we say $T$A exists.

    As before, we use $\text{acl}$ and $\ind$ to denote the algebraic closure and non-forking independence in $T$; while we use $\text{acl}_\sigma$ and $\indi{\sigma}$ to denote those in $T$A.
	\begin{fact}[{\cite[3.5, 3.6 and 3.8]{CHATZIDAKIS199871}}]\label{F: TA}
		Suppose $T$ is stable, complete, and admits QE and EI, and assume that $T$A exists. Then, the following holds:
		\begin{enumerate}
			\item If $(M_1,\sigma_1),\,(M_2,\sigma_2)$ are two models of $T$A, containing a common substructure $(A,\sigma)$ which is algebraically closed in the sense of $T$, then $(M_1,\sigma_1)\equiv_A(M_2,\sigma_2)$.
			\item Let $A\subset (M,\sigma)\models T$A. The algebraic closure of $A$ in $T$A is given by $$\text{acl}_\sigma(A)=\text{acl}\left(\bigcup_{i\in\mathbb{Z}}\sigma^i(A) \right).$$
			\item For algebraically closed sets, the non-forking independence relation in $T$A coincides with that in $T$. That is, suppose $A,B,C\subset(M,\sigma)\models T$A such that $A=\text{acl}_\sigma(A),\, B=\text{acl}_\sigma(B),\,C=\text{acl}_\sigma(C)$. Then $$A\ind^{\sigma}_C B\text{ if and only if } A\ind_C B.$$
			\item Every completion of $T$A is simple. And if $T$ is superstable, then every completion of $T$A is supersimple.
		\end{enumerate}
	\end{fact}
    \begin{defi}
        Let $T$ be a complete $L$-theory. We say $T$ admits geometric elimination of imaginaries (\textbf{gEI}) if for every $M\models T$ and every element $e\in M\eq$, there exists a finite tuple from $M$ which is interalgebraic with $e$.
    \end{defi}
	\begin{fact}[{\cite[Proposition 4.5 and 4.7]{gpoid}}]\label{F:TAGEI}       
        For $T=T\eq$ superstable such that $T$A exists, $T$A admits gEI.
    \end{fact}
        Hrushovski \cite{gpoid} proves an equivalent condition, namely 3-uniqueness, for every completion of $T$A to have elimination of imaginaries for $T=T\eq$. Following his method (see also {\cite[(4) on page 4488]{bays2017model}}), one can show Fact 1.5 without further assumptions on $T$.
        
        It is well-known that if $T$A exists then so does $T\eq$A, and $T\eq $A$\subset(T\text{A})\eq$ naturally. We can apply Fact 1.5 to $T\eq$ to obtain that $T\eq$A admits gEI. Thus, as any imaginary $e$ in a model of $T$A can be viewed as an imaginary in $T\eq$A, it is interalgebraic with a tuple from $T\eq$. If $T$ has gEI in some sorts $\mathcal{S}$, there is a finite tuple from sort $\mathcal{S}$ which is interalgebraic with $e$. In particular, we may infer the following: 
    \begin{rmk}\label{R:TAgEI}   
        Let $T$ be a complete superstable theory such that $T$A exists. Suppose that $T$ has gEI in some sorts $\mathcal{S}$. Then, $T$A also has gEI in $\mathcal{S}$.
    \end{rmk}
    \begin{fact}[{\cite[Proposition 4.8]{nonfirstorder}}]\label{F:TAP=TPA}
        Let $T$ be a complete stable theory. If $T$A exists, then $T_P$ is first-order and $T_P$A exists. Moreover, let $T'$ be any completion of $T$A, then $T'_P$, the theory of lovely pairs of models of $T'$, is a completion of $T_P$A and first-order. In fact, every completion can be obtained in this way. In this sense, $T_P$A and $T$A$_P$ coincide.
    \end{fact}

    \begin{exmp}
        Let $T=$ACF$_p$ for $p$ prime or zero. $T$ is known to be strongly minimal, so in particular $\omega$-stable and it follows from \cite{MR704080} that $T_P$ is first-order. In fact, $T_P$ is the theory of proper pairs of algebraically closed fields. In particular, every $\aleph_1$-saturated model of $T_P$ is of the form $(K,P(K))$ where $K$ has uncountable transcendence degree over $P(K)$, and $P(K)$ has uncountable transcendence degree over the prime field. Moreover, $T_P$ is $\omega$-stable by Fact~\ref{F:TP}.3.
        
        By \cite{macintyre1997generic}, $T$A exists (also known as ACFA). By Fact~\ref{F: TA}.4, any completion of ACFA is supersimple and admits gEI by Remark~\ref{R:TAgEI}. Actually, every completion of ACFA admits full elimination of imaginaries, see \cite[1.10]{ZHdiffield}. By Fact~\ref{F:TAP=TPA}, $T_P$A exists, and $T_P$A coincides with $T$A$_P$, the theory of lovely pairs of (completions of) ACFA. In particular, ACFA$_P$ is first-order.
    \end{exmp}
	\subsection*{Pseudofinite fields}
	\begin{defi}
		Let $Q$ denote the set of all prime powers. For $q\in Q$, let $\mathbb F_q$ denote the finite field with $q$ elements. We define the $L_\text{ring}$-theory
		$$\text{Psf}:=\text{Th}(\{\prod_\mathcal{U}\mathbb{F}_q\mid \mathcal{U}\text{ is a non-principal ultrafilter on } Q\})$$as the common theory of non-principal ultraproducts of finite fields. Its models are called \textbf{pseudofinite fields}.
	\end{defi}
	In \cite{ax1968elementary}, Ax gave a characterization of Psf. A field $F$ is pseudofinite if and only if it satisfies the following conditions:
	\begin{enumerate}
		\item Gal$(F)\cong\hat{\mathbb{Z}}$.
		\item $F$ is perfect.
		\item $F$ is pseudo-algebraically closed, i.e., every (non-empty) absolutely irreducible variety defined over $F$ has an $F$-rational point.
	\end{enumerate}
	These properties can be expressed as a first-order theory, providing an axiomatization of Psf.\\
	
	Every pseudofinite field $F$ can be embedded into a model $(K,\sigma)$ of ACFA as the fixed field of $K$ \cite{afshordel2014generic}. Since any completion of ACFA is (super-)simple, it follows that any completion of Psf is also (super-)simple. Since ACFA admits EI, $F$ is stably embedded in $(K,\sigma)$, i.e., every definable subset of $F$ in $(K,\sigma)$ is definable using parameters from $F$ \cite{ZHdiffield}. (Actually, it is definable in the pure field language by Fact~\ref{F:indstronF} below.) However, $F$ still retains extra structure from $\sigma$.

    For instance, the $2n$-ary relation $R_n(\bar x,\bar y)$ which asserts that the polynomial $T^n+x_1T^{n-1}+\dots+x_n$ is irreducible over $F$, and that for any root of the polynomial $\alpha$, $\sigma(\alpha)=y_1+y_2\alpha+\dots+y_n\alpha^{n-1}$. It's easy to see that $R_n(\bar x,\bar y)$ is $\emptyset$-definable in $L_\sigma$, but only definable in $L$ with parameters for $n>4$. Fortunately, adding $\sigma$ does not increase much complexity on $F$ as Fact 1.9 below shows that $\code{R_n}$'s are all the extra information obtained from $\sigma$.\footnote{For a definable set $X$, we denote by $\code{X}$ the canonical parameter of $X$.}
    \begin{nota}
        For a pseudofinite field $F$ and $n>1$, there is a unique field extension of degree $n$ over $F$ up to isomorphism. We denoted it by $F_n$.
    \end{nota}
    Notice that if $F$ is the fixed field of $(K,\sigma)\models$ACFA, the field $E_n:=\mathrm{Fix}(\sigma^n)$ has degree $n$ over $F$, i.e., $F_n\cong_F E_n$.
    \begin{fact}[{\cite[Proposition A. page 3014]{ZHdiffield}}]\label{F:indstronF}
        The induced structure on $F$ is precisely the field structure together
    with the distinguished (imaginary) constants $\{\code{R_n}\}_{n\in\mathbb{N}^*}$, where $\code{R_n}$ codes the isomorphism type of $(F_n,\sigma)$. Moreover, as $\text{Gal}(F_n/F)$ is finite, $\code{R_n}\in\aeq{\emptyset}$ for every $n$.
    \end{fact}

    Moreover, let $A,B,C\subset F\models$Psf, which is the fixed field of some $(K,\sigma)\models$ACFA. By Fact~\ref{F: TA}.3, as $A,B,C$ are already closed under $\sigma$, the non-forking independence relation of $A,B,C$ in $K$ is given by algebraic independence. And by Fact~\ref{F:indstronF}, since the extra structure given by $\sigma$ is fully contained in $\aeq{\emptyset}$, it has no effect on the non-forking independence in $F$. We obtain the following:

    \begin{fact}\label{F:forkinginF}
        In pseudofinite fields, the non-forking independence relation is given by algebraic independence. In particular, Psf is supersimple of SU-rank 1.
    \end{fact}

    In \cite{VASSILIEV2003103}, it is shown that if a theory $T$ is simple of SU-rank 1, then $T_P$ is first-order. It follows:
    \begin{fact}
        The theory of lovely pairs of Psf is first-order. 
    \end{fact}	
	\subsection*{Pillay’s geometric description of imaginaries in pairs of ACF}
	In 2007, Pillay gave a geometric description of imaginaries in lovely (beautiful) pairs of ACF. He showed that every imaginary in lovely pairs of ACF is interalgebraic with some imaginary of a certain algebraic form \cite{PILLAY200713}.
    \begin{defi}
		Let $G$ be a connected algebraic group and $V$ be an absolutely irreducible variety, both defined over some algebraically closed field $k$. A rational function $\mu:G\times V\dashrightarrow V$ defined over $k$ is called a \textbf{generically free rational group action} if for any $K\supset k$, any $v\in V(K)$ generic (i.e., of maximal transcendence degree) over $k$, and any $g_1,g_2\in G(k)$, we have $(g_1,v),(g_2,\mu(g_1,v))\in \text{dom} (\mu)$, $\mu(g_1,\mu(g_2,v))=\mu(g_1*g_2,v)$ and $\mu(g_1,v)=v$ if and only if $g_1=\text{id}$.
	\end{defi}
    \begin{defi}
		Let $T$ be a complete expansion of an $L_{\text{ring}}$-theory of perfect PAC fields. Let $(M,P)$ be a proper elementary pair of $T$. 
        
        Let $B\subset\A^n, V\subset \A^m$ be irreducible varieties over the prime field together with a surjective morphism $\pi:V\rightarrow B$ such that for any $b\in B$ generic over the prime field, $V_b:=\pi^{-1}(b)$ is an absolutely irreducible variety defined canonically over $b$. Let $G\rightarrow B$ be an algebraic group over $B$, i.e., every fiber $G_b$ is an algebraic group and the group operations are defined uniformly over $b$.
        
        Suppose moreover there exists an $\emptyset$-definable subset $\mu\subset G\times_B V\times_B V$ such that for every $b\in B(P)$, the fiber $\mu \cap G_b\times V_b\times V_b$ is the graph of a partial map $\mu_b: G_b\times V_b\dashrightarrow V_b$, and that for generic $b\in B(P)$, $\mu_b$ is a generically free rational group action.
		
		Let $\phi(x,y)$ be the $L_P$-formula expressing the following: $y\in B(P),x\in V_y$ and for any $g_1,g_2\in G_y(P)$,
		$$(g_1,x),(g_2,\mu_y(g_1,x)) \in \text{dom} (\mu_y),$$  $$\mu_y(g_1,\mu_y(g_2,x))=\mu_y(g_1*g_2,x)\in V_y,$$
		$$\mu_y(g_1,x)=x\text{ if and only if }g_1=\text{id}_{G_y}.$$
		Set
		$$W:=\{x\in V\mid M\models\phi(x,\pi(x))\}.$$
		
		We define a binary relation $\sim$ on $W$ by setting $w_1 \sim w_2$ if $\pi(w_1)=\pi(w_2)$ and there exists $g\in G_{\pi(w_1)}(P)$ such that $\mu (g,w_1)=w_2$. One can check that $\sim$ defines an equivalence relation on $W$. We denote the quotient sort by
		$$\mathcal{S}(B,V,G,\mu):=W/\sim.$$
		
		An imaginary $e\in(M,P)\eq$ is said to be of \textbf{Pillay form} if it is a ``generic element'' in some Pillay sort $\mathcal{S}(B,V,G,\mu)$. That is, there exists a generic $b$ in $B(P)$ and $a\in V_b(M)$ a generic over $P$, such that $e=\lceil G_b(P)*a\rceil$. Here, generic means Loc$(b)=B$, Loc$(a/P)=V_b$. We denote by $(K,P)^\mathcal{P}$ the imaginaries of Pillay form in $(K,P)$.
	\end{defi}
    \begin{fact}[\cite{PILLAY200713}]\label{F:PillaygEI}
		Let $(K,P)\models($ACF$_p)_P$, where $p$ is either a prime or zero. For any imaginary $e\in(K,P)\eq$, there exists an imaginary $e'$ of Pillay form such that $e$ and $e'$ are interalgebraic.
	\end{fact}
	\begin{rmk*}       
        In the original statement of the theorem, it is not required that $b$ is the canonical parameter of $V_b$. However, from the proof in \cite{PILLAY200713}, we can assume it.
	\end{rmk*}
	Fact~\ref{F:PillaygEI} implies that (ACF$_p)_P$ admits gEI within the Pillay sorts. In this article, we will show that, after adding countably many constants, every imaginary in Psf$_P$ is interalgebraic with an imaginary of Pillay form in some finite extensions. Thus, the theory of lovely pairs of pseudofinite fields also admits gEI in Pillay sorts in finite extensions.
    
    In \cite{BoxallThesis}, it is shown that in Fact~\ref{F:PillaygEI}, one can choose $e'$ such that $e\in\deqP{e'}$, i.e., we have weakly elimination of imaginaries down to Pillay sorts.

	\section{Imaginaries in Psf$_P$}
	Let $(F,P(F))$ be a lovely pair of Psf. We will prove an analogue of Pillay's geometric description on imaginaries for $(F,P(F))$ within the framework of a richer structure $(K,P,\sigma)$, a lovely pair of ACFA. We begin by showing that $(F,P(F))$ can be embedded nicely into $(K,P,\sigma)$.
    
    As mentioned above, both $\mathrm{Psf}_P$ and $\mathrm{ACFA}_P$ are first-order, and every pseudofinite field $F$ can be realized as a stably embedded fixed field in a model of ACFA. In fact, the same observation continues to hold when one passes to lovely pairs of Psf and ACFA.
	\begin{prop}\label{P:FPinKP}
		For $(K,\sigma)\models$ACFA, suppose $(K,\sigma,P)$ is a lovely pair of Th$(K,\sigma)$. Let $F=K^\sigma$ denote the fixed field, which itself is a model of Psf. Then, $(F,P(F))$ is also a lovely pair of Th$(F)$.
	\end{prop}
	\begin{proof}
		We verify the definition directly:
		
		First, since $P(K)\prec K$ in the sense of $L_\sigma$, the corresponding fixed fields then satisfy $P(F) \prec F$ in $L_\sigma$, and hence also in $L$.
		
		For any $L$-type $p\in S_x(A)$ for some small $A\subset F$, let $p'\in S_x(A)$ be an $L_\sigma$-extension of $p$ and containing the formula $\sigma(x)=x$ (which is consistent with $p$). Since $(K,\sigma,P)$ is a lovely pair, there exists a non-forking extension $q'\in S(A\cup P(K))$ of $p'$ which is realized in $K$, and hence in $F$, since $\sigma(x)=x\,\in q'$. Let $q$ be the $L$-reduct of $q'$ over $A \cup P(F)$. Thus, $q$ is a non-forking extension of $p$ realized in $F$ as non-forking in Psf coincides with non-forking in $K$ by Fact~\ref{F:forkinginF}. Hence, the heir property is satisfied.
		
		Now we verify the coheir property. Again, let $p\in S_x(A)$ be an $L$-type as before. If $p$ does not fork over $P(A)$, take any $p'\in S_x(A)$ an $L_\sigma$-extension of $p$ such that $\sigma(x)=x\in p'$. Since $A,P(A)$ are both fixed by $\sigma$ and $p'$ concentrates on $F$, the fixed field of $\sigma$, it follows from Fact~\ref{F: TA}.3 and Fact~\ref{F:forkinginF} that $p'$ does not fork over $P(A)$. Then, $p'$ is realized in $P(F)$, thus so is $p$. Hence, the coheir property holds. 
        Therefore, $(F,P(F))$ is a lovely pair.
	\end{proof}
	
	Let $F\models$Psf and $(K,\sigma)\models$ACFA, with $F$ being the fixed field of $K$. We take $(K',\sigma,P')$ to be a lovely pair of Th$(K,\sigma)$ and let $F'$ denote the fixed field of $K'$. It follows that $F'\equiv F$ since $K'\equiv_\sigma K$. Moreover, by Proposition~\ref{P:FPinKP}, $(F',P'(F'))$ is an $L$-lovely pair. Henceforth, we may assume that $(F,P(F))$ is a lovely pair of pseudofinite fields and $(K,\sigma,P)$ is a lovely pair of ACFA, with $F$ being the fixed field of $K$.
	\begin{prop}\label{P:GEIinK}
		$(K,P,\sigma)$ admits gEI in Pillay sorts. More precisely, for every $e\in (K,P,\sigma)^\text{eq}$, there exists an element $e'$ from $(K,P)^\mathcal{P}$ such that $e$ and $e'$ are $L_{P,\sigma}$-interalgebraic in $(K,P,\sigma)^\text{eq}$.
	\end{prop}
	\begin{proof}
		By Fact~\ref{F:TAP=TPA}, we know $(K,P,\sigma)$ is a model of (ACF$_P$)A. By Fact~\ref{F:PillaygEI}, $(K,P)$ has gEI in $(K,P)^\mathcal{P}$. Then, applying Remark~\ref{R:TAgEI}, we obtain the desired statement.
	\end{proof}
	For $e\in (F,P)^\text{eq}$, we can regard $e$ as an element in $(K,P,\sigma)\eq$. Then, from Proposition~\ref{P:GEIinK}, we are able to find $e'=\lceil G_b(P)*a\rceil\in(K,P)^\mathcal{P}$ of Pillay form such that $e'$ is $L_{P,\sigma}$-interalgebraic with $e$ for some $a,b\in K$. Moreover, the parameters $a$ and $b$ essentially come from the fixed field $F$ since $e'$ is algebraic over $e$, and $e$ lies in the definable closure of the fixed field. Let $F\alg$ denote the (field-theoretic) algebraic closure of $F$. We will show that $e'\in\deqP{F\alg}\cap(K,P)^\mathcal{P}$.
	\begin{prop}\label{P:Falgcls}
		Let $(F\alg,P)\eq:=\deqP{F\alg}$ denote the definable closure of $F\alg$ in $(K,P)\eq$ in the sense of the language $L_P\eq$. Then, $\sigma$ is still an automorphism of $(F\alg,P)\eq$. And as an $L\eq_{P}\cup\{\sigma\}$-substructure, $(F\alg,P)\eq\subset(K,P)\eq$ is algebraically closed.
	\end{prop}
	\begin{proof}
		By Fact~\ref{F: TA}.2, it suffices to show that $(F\alg,P)\eq$ is algebraically closed in $(K,P)\eq$ in the sense of $L_P\eq$ and closed under $\sigma$.
		
		For any $a\in F\alg$ there is $f\in F[T]$ such that $f(a)=0$, hence $f(\sigma(a))=0$. Hence $\sigma(a)\in F^{\text{alg}}$; that is, $F^{\text{alg}}$ is closed under $\sigma$. Moreover, since the predicate $P$ is fixed setwise by $\sigma$, $(F\alg,P)\eq$ is also closed under $\sigma$.
		
		By the definition of lovely pairs, $P(K)$ is an elementary submodel of $K$ in the sense of ACFA. For $a_1\dots a_n\in P(K)$, if there exist $b_1,\dots, b_n\in E_N$ for some $N$ such that $\sum_ia_ib_i=0$, then by elementarity, there also exist $b_1'\dots b_n'\in P(E_N)$ such that $\sum_ia_ib_i'=0$. We conclude that $P(K)$ is linearly disjoint from $E_N$ over $P(E_N)\subset P(F\alg)$. And since $F\alg=\bigcup_{N\in\mathbb{N}}E_N$, $F\alg\ind_{P(F\alg)}P(K)$ holds in the sense of fields. Consequently, for any tuple $d\in F\alg$, $\Cb{\tp{d/P(F\alg)}}=\Cb{\tp{d/P(K)}}$. It follows from Fact~\ref{F:TP}.1 that $(F\alg,P)\prec(K,P)$ in the sense of $L_P$. Therefore, $(F\alg,P(F\alg))\eq$ is algebraically closed in $(K,P)\eq$ in the sense of $L_P\eq$ and hence in the sense of $L_P\eq\cup\{\sigma\}$.
	\end{proof}
	
	As $e\in(F,P(F))\eq$, choose $a_0\in F$ such that $e\in\text{dcl}_{P,\sigma}\eq (a_0)$ in $(K,P,\sigma)$. We now apply Proposition~\ref{P:Falgcls} to the Pillay form imaginary $e'$ obtained above.
    
    We can see $a_0$ and $e'$ as two real elements in the language $L':=L_{P}\eq\cup\{\sigma\}$. Since $e'\in\aeqs{a_0}$ in $L_{P,\sigma}\eq=(L')\eq$, it follows $e'$ is algebraic over $a_0$ in $L'$. By Proposition \ref{P:Falgcls}, $\deqP{F\alg}$ is $L'$-algebraically closed and $e'\in\aeqs{a_0}$, we conclude $e'\in\deqP{F\alg}$. Since $e'$ is of Pillay form, we may choose its parameter $a,b$ from $F\alg$. As $F\alg=\bigcup_{n\in \mathbb N} E_n$ there exists $N\in\mathbb{N}$ such that $a,b\in E_N$. The following lemma shows that we can understand $e'$ as a Pillay form element in $(E_N,P(E_N))^\mathcal{P}$.
	\begin{lemma}\label{L:FalgtoFN}
		Let $e'=\lceil G_b(P)*a\rceil\in(F\alg,P)^\mathcal{P}$ be of Pillay form such that $a,b\in E_N$. As $E_N$ is the fixed field of $\sigma^N$, the orbit $G_b(P(E_N))*a$ is definable in $(K,P,\sigma)$. Denote its code by $e''$. Then, $e'$ and $e''$ are interdefinable in the full language.
	\end{lemma}
	\begin{proof}
		Suppose $e'\in\mathcal{S}(B,V,G,*)$ in $(F\alg,P(F)\alg)\prec(K,P)$. Let $X$ (not necessarily definable) denote the set of all imaginaries of the form $\lceil G_y(P(E_N))*x\rceil$ for some $y\in B(P(E_N))$ and some $x\in V_y(E_N)$ generic (i.e., of maximal transcendence degree) over $P(F\alg)$. There is a natural $\emptyset$-$L_\sigma$-definable quotient map from $X$ to $\mathcal{S}(B,V,G,*)$, which takes $\lceil G_y(P(E_N))*x\rceil$ to $\lceil G_y(P)*x\rceil$. So it takes $e''$ to $e'$. It suffices to show that this map is injective.
		
		Suppose $x'\in V_y(E_N)$ generic over $P$ and $g\in G_y(P)$ are such that $g*x'=x$. Since $(F\alg,P(F)\alg)\prec(K,P)$ and the group action is generically free, $g\in G_y(P(F\alg))$. For every $\tau\in\text{Gal}(F\alg/E_N)$, we have $\tau(g)*x'=x$, since $x,x'\in E_N$ and the group action is defined over $y\in E_N$. As the group action is generically free, we have $g=\tau(g)$. Hence $g\in G_b(P(E_N))$. That is to say $\lceil G_b(P(E_N))*x\rceil=\lceil G_b(P(E_N))*x'\rceil$.
	\end{proof}
    We have thus reduced the original imaginary in $(F,P)^{\mathrm{eq}}$ to an imaginary of Pillay form in a finite extension $(E_N,P(E_N))$, with interalgebraicity computed in the ambient structure $(K,P,\sigma)$. To bring this description back to the original pair $(F,P)$, it remains to realize these finite extensions uniformly inside $(F,P)$. We do this by interpreting, after naming suitable parameters, a copy of each $E_N$ in $F$.

    
\subsection*{Interpretation of finite extensions}

	Let $C=\{c_{i,j}\}_{i\geq 2,1\leq j\leq i}\subseteq P(F)$ be such that $f_N(T):=T^N+\sum_{j=1}^N c_{N,j}T^{N-j}$ is an irreducible polynomial of degree $N$ over $P(F)$ for all $N>1$. Hence $F_N\cong F[T]/f_N(T)$, $F_N$ is interpretable in $F$ over $C$. With this interpretation, every element in $(F_N,P(F_N))^\mathcal{P}$ can be viewed as an element in $(F,P,C)\eq$.
    
    For every $N>1$, fix $t_N\in P(F\alg)$ a root of $f_N$. We have a bijection $\iota_{t_N}:\:F_N\rightarrow E_N$ defined by sending $\bar T$ to $t_N$. Each $\iota_{t_N}$ gives an embedding of $F_N$ into $F\alg$. Notice that these embeddings preserve the predicate $P$, and $\iota_{t_N}$ also induces an injection from $(F_N,P(F_N))\eq\subset(F,P,C)\eq$ to $(K,P,\sigma)\eq$ by sending an equivalence class in $(F_N,P(F_N))$ to an equivalence class in $(K,P,\sigma)$. 

	\begin{thm}\label{T:FhasGEIinK}
		Let $(F,P)\models\mathrm{Psf}_P$. Let $(K,\sigma,P)$ be a model of ACFA$_P$ such that $F=$Fix$_\sigma (K)$. Then, with the notation above, for any $e\in(F,P)\eq$, there exist $N>1$ and $\tilde e\in(F_N,P(F_N))^\mathcal{P}\subset(F,P,C)\eq$ which is interalgebraic with $e$ (in $(K,P,\sigma)\eq$) over $C$.
	\end{thm}
	\begin{proof}
		In what follows, interalgebraicity and interdefinability are always understood in the ambient structure $(K,P,\sigma)\eq$.
        
        Let $e\in(F,P)\eq$. By the previous discussion, there exist some $a,b\in E_N$ for some integer $N>1$ such that $e'':=\code{G_b(P(E_N))*a}$ is of Pillay form and $e$ is interalgebraic with $e''$. Take $\tilde e$ to be $\iota_{t_N}^{-1}(e'')\in(F_N,P(F_N))\eq\subset(F,P,C)\eq$. 
        
        As $\iota_{t_N}$ is definable over $t_N$ and there are only finitely many roots of $f_N$, we conclude that $\tilde e$ and $e''$ are interalgebraic over $C$. Therefore, $\tilde{e}$ is interalgebraic with $e'$ over $C$, hence also with $e$.
	\end{proof}
	Thus, after naming $C$, every imaginary of $(F,P)$ is geometrically represented by a Pillay form in a finite extension. However, the interalgebraicity is in the sense of the full language $L_{P,\sigma}$ in $K$, rather than in the sense of Psf$_P$. We will study the structure induced on $(F,P)$ from $(K,P,\sigma)$ and address this issue in the next section.

	\section{The induced structure on $(F,P)^{\mathrm{eq}}$}
    Let $(K,\sigma)\models$ACFA and $F$ be its fixed field. Recall that in Fact~\ref{F:indstronF}, it is stated that the induced structure of $(K,\sigma)$ on $F$ is given by only naming $\code{R_n}\in\aeq{\emptyset}$.

    Let $F_{\text{ind}}=(F,\code{R_n})_{n\in\mathbb{N}}$ denote the induced structure of $F\subset(K,\sigma)$ (by adding $R_n$ to the language). As $F$ has SU-rank 1 in ACFA and $\code{R_n}\in\aeq{\emptyset}$, Psf$_{\text{ind}}$ is still supersimple of SU-rank 1 after adding codes for definable sets. By \cite{VASSILIEV2003103}, (Psf$_\text{ind})_P$ is first-order.
	\begin{prop}\label{P:indstronFP}
		 Let $(F,P)_\text{ind}$ denote the induced structure on $(F,P(F))\subset(K,P,\sigma)$. Then, $(F_{\text{ind}},P)=(F,P)_\text{ind}$.
	\end{prop}
	\begin{proof}
		For convenience of notation, we regard $P$ as a predicate symbol rather than a subfield. We may assume $(F,P)$ is sufficiently saturated and strongly homogeneous.
        
		Notice that $(F_\text{ind},P)$ is definable in $(F,P)_\text{ind}$ since $(F_\text{ind},P)$ is definable in $(K,\sigma,P)$. To show $(F,P)_\text{ind}$ is also definable in $(F_\text{ind},P)$, we will prove that: for any tuples $a,a'\in F$ such that $a$ and $a'$ have the same type in $(F_\text{ind},P)$, they also have the same type in $(F,P)_\text{ind}$.
		
		Since $a$ and $a'$ have the same type in $(F_\text{ind},P)$, there exists $\tau\in$ Aut$_P(F_{\text{ind}})$, such that $\tau(a)=a'$. We can extend $\tau$ to $\tilde\tau\in$ Aut$_{P}(F\alg)$, and any such extension $\tilde\tau$ will commute with $\sigma$: for any $x\in F\alg$, let $T^n+\sum_{i=1}^nd_{i}T^{n-i}$ be the minimal polynomial of $x$ over $F$. Then, there exists a unique tuple $\bar b\in F$ such that $F\models R_n(\bar d,\bar b)$. i.e. $\sigma(x)=\sum_{i=1}^n b_ix^{i-1}$. Hence, $\tilde\tau(\sigma(x))=\sum_{i=1}^n \tau(b_i)\tilde\tau(x)^{i-1}$. Notice that $\tau$ is an automorphism of $F_\text{ind}$. We have $F\models R_n(\tau(\bar d ),\tau(\bar b))$. And $\tau(\bar d)$ are the coefficients of the minimal polynomial of $\tilde\tau(x)$ over $F$, therefore, $\tilde\tau(\sigma(x))=\sum_{i=1}^n \tau(b_i)\tilde\tau(x)^{i-1}=\sigma(\tilde\tau(x))$.
		
		From the above, $\tilde\tau\in$Aut$_{P,\sigma}(F\alg)$. Since $F\alg$ is algebraically closed in the sense of $L_{P,\sigma}$ and $\deqP{F\alg}$ is $L_P\eq$-algebraically closed as $(F\alg,P(F\alg))$ is a model, we can apply Fact~\ref{F: TA}.1 to conclude that $\tilde\tau$ extends to Aut$_{P,\sigma}(K)$. So, $a$ and $a'$ have the same type in $(K,P,\sigma)$ hence also in $(F,P)_\text{ind}$.
	\end{proof}
    
	\begin{coro}\label{C:algclscoincd}
		For $e,e'\in(F,P)\eq$, if $e$ and $e'$ are interalgebraic in the sense of $L_{P,\sigma}$ in $(K,P,\sigma)\eq$, then $e$ and $e'$ are interalgebraic in the sense of $L_P$ in $(F,P)\eq$.
	\end{coro}
	\begin{proof}
		By Fact~\ref{F:indstronF}, $F_\text{ind}$ is obtained by adding all the codes $\{\lceil R_n\rceil\}_{n\in\mathbb N^*}$. From Proposition~\ref{P:indstronFP}, we know that the induced structure of $(F,P)$ can also be obtained from $\lceil R_n\rceil$'s. So every automorphism of $(F,P)$ fixing all $\lceil R_n\rceil$ is also an automorphism of $(F,P)_\text{ind}$. Hence, for any imaginary $e\in(F,P)\eq$, the algebraic closure of $e$ in $(F,P(F))\eq$ in the full language is given by: 
        $$\aeqs{e}\cap(F,P(F))\eq=\aeqP{e\cup\bigcup_n\{\lceil R_n\rceil\}}\cap(F,P(F))\eq.$$
       Since $\lceil R_n\rceil\in\aeq{\emptyset}$, we have $\aeqs{e}\cap(F,P(F))\eq=\aeqP{e}\cap(F,P(F))\eq$, (the algebraic closures are calculated in Psf$_P$). Thus if two imaginaries from $(F,P)\eq$ are interalgebraic in $(K,P,\sigma)\eq$, then they are interalgebraic in $(F,P)\eq$.
	\end{proof}

    As a complementary, we show that $(F,P)$ is stably embedded in $(K,P,\sigma)$, even though we do not need it for the main results of this paper.
    \begin{prop}\label{P:stblyembed}
        Let $(K_0,P(K_0),\sigma)\preceq(K,P,\sigma)$ be a small elementary submodel. Denote by $F_0$ the fixed field of $K_0$. If $a,a'\in F$ have the same $L_{P,\sigma}$-type over $F_0$, then they have the same type over $K_0$. In other words, $(F,P(F))$ is stably embedded in $(K,P,\sigma)$.
    \end{prop}
    \begin{proof}
    Notice that since $F$ is linearly disjoint from $P$ over $P(F)$, we have $F$ is linearly disjoint from $F_0P$ over $F_0P(F)$, thus $\Cb{\mathrm{tp}_L(a/F_0P)}\subset F_0P(F)\subset F$. Replacing $a$ by $(a,\Cb{\mathrm{tp}_L(a/F_0P)})$, we may assume that $aF_0$ is $P$-independent.
    
    Since $K_0\preceq K$ (as $L_{P,\sigma}$-structures), $\tP{F/K_0}$ is finitely satisfiable in $F_0$. So we have $K_0\indi{P}_{F_0}F$ and hence $K_0\ind_{F_0P}a$. As $K_0\ind_{P(K_0)}P$, $K_0\ind_{F_0P(K_0)}F_0P$. We conclude $K_0\ind_{F_0P(K_0)}aF_0P$. So, $K_0a\ind_{aF_0P(K_0)}P$. Notice that $aF_0\ind_{P(aF_0)}P$, thus by transitivity, $K_0a\ind_{P(K_0a)}P$. Therefore, $(K_0a)\alg$ is $L_P$-algebraically closed.
    
    As $F$ is stably embedded in $(K,\sigma)$ (see \cite[Proposition 1.11]{ZHdiffield}), there exists $\tau\in\mathrm{Aut}_{\sigma}(K)$ such that $\tau$ fixes $K_0$ and sends $a$ to $a'$. Let $\tau_0$ be the restriction of $\tau$ on $(K_0a)\alg$. Then, $\tau_0$ is an $L_{P,\sigma}$-isomorphism from $(K_0a)\alg$ to $(K_0a')\alg$. Since $(K_0a)\alg$ is $L_P$-algebraically closed, $((K_0a)\alg,P(K_0a)\alg)$ is an $L_P$-elementary submodel of $(K,P(K))$. Hence $\deqP{(K_0a)\alg}$ is $L_P\eq$-algebraically closed. Therefore, by Fact~\ref{F: TA}.1, $a$ and $a'$ have the same $L_{P,\sigma}$-type over $K_0$.

    \end{proof}

	\begin{thm}[Main Theorem]
		Let $(F,P)$ be a model of lovely pairs of pseudofinite fields. We take a parameter set $C=\{c_{ij}\}_{0<j\leq i,i>1}$ such that $f_N(T):=T^N+\sum_{j=1}^Nc_{Nj}T^{N-j}$ is an irreducible polynomial of degree $N$ over $P(F)$ for all $N>1$. Then, for any $e\in(F,P)\eq$, there exist $N>1$ and $\tilde e\in(F_N,P(F_N))^\mathcal{P}\subset(F,P,C)\eq$ interalgebraic with $e$ over $C$.
	\end{thm}
	\begin{proof}    
		Theorem~\ref{T:FhasGEIinK} asserts that, over $C$, any $e\in(F,P)\eq$ is interalgebraic with some $\tilde e\in(F_N,P(F_N))^\mathcal{P}$ in the ambient structure $(K,P,\sigma)\eq$. However, by Corollary~\ref{C:algclscoincd}, the algebraic closure in the sense of the pairs of pure fields structure $(F,P)$ is the same as the algebraic closure in the induced structure $(F,P)_\text{ind}$. Therefore, $e$ and $\tilde{e}$ are interalgebraic over $C$ in $(F,P)$.
	\end{proof}
    The proof of the Main Theorem is now complete. We conclude this section with one further consequence of the comparison above: algebraic closure and non-forking in $\mathrm{Psf}_P^{\mathrm{eq}}$ can be related to their counterparts in $\mathrm{ACF}_P^{\mathrm{eq}}$.

    For any $e\in(F,P(F),C)\eq$, by the Main theorem, choose $e'\in(F_N,P(F_N))^\mathcal{P}$ for some $N$ interalgebraic with $e$ (over $C$). For any root $t_N$ of $f_N$, $\iota_{t_N}(e')\in\text{dcl}\eq_{P,\sigma}(F\alg)\subset(K,P,\sigma)\eq$. Moreover, by Lemma~\ref{L:FalgtoFN}, there exists an imaginary in $\deqP{F\alg}\subset(K,P,\sigma)\eq$, denoted by $\nu(e)$, such that $\iota_{t_N}(e')$ is $L_{P,\sigma}\eq$-interdefinable with $\nu(e)$. Notice that $(F\alg,P(F\alg))\prec(K,P)$ as an $L_P$-substructure, $\nu(e)$ is an element in $(F\alg,P(F\alg))\eq$ and $\nu(e)$ is interalgebraic with $e$ in $(K,P,\sigma)\eq$.

    For each $e\in(F,P(F),C)\eq$, choose once and for all such an imaginary $\nu(e)$. Although $\nu(e)$ depends on these choices, its algebraic-closure class over $C$ does not, which is all that will be used below.
    \begin{prop}\label{P:reductprncpl}
        Let $(F,P(F))\models$Psf$_P$, $(F\alg,P(F)\alg)\models$ACF$_P$. Let $C$ and $\nu$ be as before. Let $A,B,D\subset (F,P)\eq$ be collections of imaginaries. Then,
        \begin{enumerate}
            \item $B\subset\mathrm{acl}_{(F,P(F),C)\eq}(A)\text{ if and only if }\nu(B)\subset\mathrm{acl}_{(F\alg,P(F\alg),C)\eq}(\cup_{i\in\mathbb{Z}}\sigma^i(\nu(A))).$
            \item $A\indi{}_{CD} B\text{ in Psf$_P\eq$ iff } \cup_{i\in\mathbb{Z}}\sigma^i(\nu(A))\ind_{C\cup_{i\in\mathbb{Z}}\sigma^i(\nu(D))} \cup_{i\in\mathbb{Z}}\sigma^i(\nu(B))$ in ACF$_P\eq$. 
        \end{enumerate}
    \end{prop}
    \begin{proof}
    Take $(K,P,\sigma)\models$ACFA$_P$ such that $(F,P(F))$ is a substructure and is the fixed field of $K$. Working over $C$, we have:
    \begin{enumerate}
        \item By Corollary~\ref{C:algclscoincd}, 
    $$\mathrm{acl}_{(F,P(F))\eq}(A)=\mathrm{acl}_{(K,P,\sigma)\eq}(A)\cap (F,P(F))\eq.$$
    Since $A$ and $\nu(A)$ are interalgebraic in the ambient structure over $C$,
    $$\mathrm{acl}_{(F,P(F),C)\eq}(A)=\mathrm{acl}_{(K,P,\sigma.C)\eq}(\nu(A))\cap (F,P(F))\eq.$$
    By Fact~\ref{F: TA}.2,
    $$\quad\mathrm{acl}_{(K,P,\sigma)\eq}(\nu(A))\cap (K,P)\eq=\mathrm{acl}_{((K,P)\eq,\sigma)}(\nu(A))=\aeqP{\cup_{i\in\mathbb{Z}}\sigma^i(\nu(A))}.$$
    Hence, $B\subset\mathrm{acl}_{(F,P(F),C)\eq}(A)$ if and only if $\nu(B)\subset\mathrm{acl}_{(K,P,\sigma,C)\eq}(\nu(A))$ if and only if $\nu(B)\subset\mathrm{acl}_{(F\alg,P(F\alg),C)\eq}(\cup_{i\in\mathbb{Z}}\sigma^i(\nu(A))).$
    \item Similarly, $A\ind_D B$ in Psf$_P\eq$ if and only if $A\ind_D B$ in (Psf$_P$)$_\text{ind}\eq$ since they only differ by some constants in $\aeq{\emptyset}$ by Proposition~\ref{P:indstronFP}. The non-dividing independence in (Psf$_P$)$_\text{ind}$ is by definition the non-dividing independence in $(K,P,\sigma)\eq$ as an ACFA$_P\eq$ model and since both theories are simple, non-forking and non-dividing coincides, we have $A\ind_{D}B$ in (Psf$_P$)$_\text{ind}$ if and only if $A\ind_{D}B$ in ACFA$_P\eq$. As $\nu(A),\nu(B),\nu(D)$ are interalgebraic with $A,B,D$ (over $C$) in ACFA$_P\eq$ respectively. The independence relation is equivalent to $\nu(A)\ind_{\nu(D)}\nu(B)$ in ACFA$_P\eq$, and thus in ACF$_P\eq$A. By Fact~\ref{F: TA}.3, the independence relation is equivalent to $\cup_{i\in\mathbb{Z}}\sigma^i(\nu(A))\ind_{C\cup_{i\in\mathbb{Z}}\sigma^i(\nu(D))} \cup_{i\in\mathbb{Z}}\sigma^i(\nu(B))$ in $(K,P)\eq$ as a model of ACF$_P\eq$.
    \end{enumerate}
    \end{proof} 
Finally, if \(A\) is finite, then \(\nu(A)\) is definable over a real tuple from some \(E_N\), and hence $\cup_{i\in\mathbb Z}\sigma^i(\nu(A))$ is finite. Thus, for finite sets of imaginaries, the comparison above involves only finitely many imaginaries in ACF$_P\eq$.

	\section{Further Remarks}
	In this paper, we showed that every imaginary in lovely pairs of pseudofinite fields is interalgebraic with an imaginary of Pillay form in the corresponding pair of finite extensions. However, our interpretation requires parameters defining an irreducible polynomial of each degree over $P(F)$. If these constants can all be chosen from the prime field, no additional constants are needed for our construction. In particular, working over a model of Psf ensures that sufficiently many constants are available.
	
	However, in general, there may not exist enough constants to define these irreducible polynomials. Indeed, there exist models of Psf containing even the algebraic closure of the prime field. It remains open whether geometric elimination of imaginaries can be established in the Pillay sorts of $(F,P)$ itself, rather than only in the finite extensions, and whether this can be done without introducing additional constants.
    
    \footnotesize
	\bibliographystyle{alpha}
	\bibliography{aty}
\end{document}